\documentclass[12pt,twoside]{article}
\usepackage[all]{xy}   
        \usepackage {amssymb,latexsym,amsthm,amsmath}
        
\long\def\symbolfootnote[#1]#2{\begingroup%
\def\thefootnote{\fnsymbol{footnote}}\footnote[#1]{#2}\endgroup}
\newcommand{\Z}{\ensuremath{\mathcal{Z}}}

\newcommand{\tr}{\ensuremath{{}^t\!}}
\newcommand{\tra}{\ensuremath{{}^t}}
\newcommand{\C}{\mathfrak C}

\newcommand{\Aut}{\textup{Aut}}
\newcommand{\diag}{\textup{diag}}
\newcommand{\gal}{\textup{Gal}}
\newcommand{\cross}{\times}
\newcommand{\tensor}{\otimes}
\makeatletter
\def\imod#1{\allowbreak\mkern10mu({\operator@font mod}\,\,#1)}
\makeatother

\newtheorem{theorem}{Theorem}[subsection]
\newtheorem{lemma}[theorem]{Lemma}
\newtheorem{corollary}[theorem]{Corollary}
\newtheorem{proposition}[theorem]{Proposition}
\newtheorem*{theorem*}{Theorem}
\theoremstyle{definition}

\numberwithin{equation}{subsection}

\newcommand{\ignore}[1]{}

\newcommand{\mynote}[1]{}
\begin{document}
\setcounter{section}{0}
\title{\bf Reality Properties of Conjugacy Classes in Algebraic Groups }
\author{Anupam Singh \hskip2mm and \hskip2mm Maneesh Thakur\\ Stat Math Unit,  Indian Statistical Institute \\
8th Mile Mysore Road, Bangalore 560059,  India\\
email : anupamk18/maneesh.thakur@gmail.com}
\date{}

\maketitle
\vskip-1.2cm
\symbolfootnote[0]{2000 Mathematics Subject Classification 20G15 (primary), 20G10, 20E45 (secondary).}
\symbolfootnote[0]{Keywords : Algebraic Groups, Real Elements, Conjugacy Classes.}

\begin{abstract}
Let $G$  be an algebraic group defined over a field $k$.
We call $g\in G$ {\bf real} if $g$ is conjugate to $g^{-1}$ and $g\in G(k)$ as {\bf $k$-real} if $g$ is real in $G(k)$.
An element $g\in G$ is {\bf strongly real} if $\exists h\in G$, $h^{2}=1$ (i.e. $h$ is an {\bf involution}) such that $hgh^{-1}=g^{-1}$.
Clearly, strongly real elements are real and are product of two involutions.
Let $G$ be a connected adjoint semisimple group over a perfect field $k$, with $-1$ in the Weyl group. 
We prove that any strongly regular $k$-real element in $G(k)$ is strongly $k$-real (i.e. is a product of two involutions in $G(k)$).
For classical groups, with some mild exceptions, over an arbitrary field $k$ of characteristic not $2$, we prove that $k$-real semisimple 
elements are strongly $k$-real.
We compute an obstruction to reality and prove some results on reality specific to fields $k$ with $cd(k)\leq 1$.
Finally, we prove that in a group $G$ of type $G_2$ over $k$, characteristic of $k$ different from $2$ and $3$, any real element in $G(k)$ is strongly $k$-real. This extends our results in \cite{st}, on reality for semisimple and unipotent real elements in groups of type $G_2$.
 
 \end{abstract}

\section{Introduction}
Let $G$ be an algebraic group defined over a field $k$.
We call an element $g\in G$ {\bf real} if $g$ is conjugate to $g^{-1}$ in $G$.
We say $g\in G(k)$ is {\bf$k$-real} if there exists $h\in G(k)$ such that $hgh^{-1}=g^{-1}$.
Note that every element in the conjugacy class of a real element $g$ is real.
Such conjugacy classes are called real.
An element $t\in G$ is called an {\bf involution} if $t^2=1$.
If an involution in $G$ conjugates $g$ to $g^{-1}$, then it follows that $g$ is a product of two involutions in $G$ and conversely, any such element is real.
An element $g\in G(k)$ is called {\bf strongly real} if $g$ is a product of two involutions in $G(k)$.

In this paper, we deal with results concerning real elements in algebraic groups, defined over an arbitrary field. 
An element $t$ in a connected algebraic group $G$ is called {\bf regular} if the centralizer of $t$ has minimal dimension (the rank of $G$), {\bf strongly regular} if its centralizer in $G$ is a maximal torus.
Let $G$ be a connected semisimple algebraic group of adjoint type defined over a perfect field $k$.
Suppose the longest element $w_0$ of the Weyl group $W(G,T)$ acts by $-1$ on the set of roots with respect to a fixed maximal torus $T$.
Then for a strongly regular element $t\in G(k)$, we prove that $t$ is real in $G(k)$ if and only if $t$ is strongly real in $G(k)$ (Theorem~\ref{streal}).
Moreover, we prove that every element of a maximal torus, containing a real strongly regular element, is strongly real. 
We show that in a split connected adjoint semisimple group $G$ defined over $k$, with $-1$ in its Weyl group, every element in a 
$k$-split maximal torus is strongly real (Proposition~\ref{splittorusreal}). 
We study the structure of real semisimple elements in  groups over fields with $cd(k)\leq 1$.
Let $k$ be such a field. 
Let $G$ be a connected reductive group defined over $k$. 
Then, semisimple elements in $G(k)$ are real in $G(k)$ (Theorem~\ref{semreal}). 
It follows that if $G$ is connected semisimple of adjoint type, with $-1$ in its Weyl group, then every semisimple element in 
$G(k)$ is strongly real in $G(k)$ (Theorem~\ref{regreal}). 
This also shows that any regular element in such a group is real. 

In later sections, we prove, with some exceptions, $k$-real semisimple elements in classical groups over a field $k$ are strongly $k$-real.
We describe these results here for convenience. For $n\not\equiv 2 \imod 4$, we prove that any $k$-real element in $SL_n(k)$ is strongly $k$-real in $SL_n(k)$ (Theorem~\ref{realsln}).
We prove that any $k$-real semisimple element in $PSp(2n,k)$ is strongly $k$-real for $n\geq 1$ (Theorem~\ref{realpsp}).
Let $Q$ be a nondegenerate quadratic form over $k$ in any dimension. 
Then $k$-real semisimple elements in $SO(Q)$ are strongly $k$-real (Theorem~\ref{realsoq}). 
Let $K$ be a quadratic extension of $k$ and let $h$ be a nondegenerate hermitian form on a $K$-vector space $V$.
We prove that $k$-real semisimple elements in $U(V,h)$ are strongly $k$-real in $U(V,h)$ (Theorem~\ref{resultsinunitary}).
We show by examples the result is false for unipotents in $U(V,h)$. 
Finally, let $G$ be a group of type $G_2$ defined over $k$, characteristic of $k$ different from $2,3$.
We prove that any $k$-real element in $G(k)$ is strongly $k$-real (Theorem~\ref{maintheorem}), this extends our results in \cite{st}. 

Our results, combined with those in \cite{pr1}, \cite{pr2}, suggest a relation between strongly real classes in groups with their orthogonal representations. 
This will be taken up in a future project  


{\bf Notation :} In what follows, we denote the centralizer of $g\in G$ by $\mathcal Z_G(g)$, the center of $G$ by $\mathcal Z(G)$ and a block diagonal matrix by $\diag(A_1,\ldots,A_n)$ where $A_i$'s are the block entries on the diagonal.
Transpose of a matrix $A$ is denoted by $\tr A$.
 

\section{Reality in Linear Algebraic Groups}
In this section we discuss reality for general linear algebraic groups.
We also compute a cohomological obstruction to reality.
We assume in this section that $k$ is a perfect field and characteristic of $k$ is not $2$.
\subsection{Strongly Regular Real Elements}

An element $t$ in a connected linear algebraic group $G$ is called {\bf regular} if its centralizer $\Z_G(t)$ has minimal dimension among all centralizers. 

An element is called {\bf strongly regular} if its centralizer in $G$ is a maximal torus.
Let $G$ be a connected, adjoint simple algebraic group defined over $k$ such that the longest element $w_0$ in the Weyl group $W$ of $G$ with respect to a maximal torus $T$ acts by $-1$ on the roots.
The adjoint groups of type $A_1, B_l, C_l, D_{2l} (l>2), E_7, E_8, F_4, G_2$ are precisely the simple groups which satisfy the above hypothesis. 
For the groups of the above type we record below a theorem of Richardson and Springer (\cite{rs}, Proposition 8.22) which plays an important role in our investigation.
\begin{theorem}[Richardson, Springer]\label{risp}
Let $G$ be a simple adjoint group over an algebraically closed field $k$. Let $T$ be a maximal torus of $G$ and let $c\in W(T)$ be an 
involution. Then there exists an involution $n\in N(T)$ which represents $c$.
\end{theorem}
We have,
\begin{theorem}\label{streal}
Let $G$ be a connected semisimple adjoint group defined over a field $k$ (not assumed algebraically closed), with $-1$ in its Weyl group.
Let $t\in G(k)$ be a strongly regular element. 
Then $t$ is real in $G(k)$ if and only if $t$ is strongly real in $G(k)$.
Moreover, every element of a maximal torus, which contains a strongly regular real element, is strongly real in $G(k)$.
\end{theorem}
\noindent{\bf Proof :} Let $t\in G(k)$ be a strongly regular real element and let $g\in G(k)$ be such that $gtg^{-1}=t^{-1}$.
Let $T$ be a maximal torus in $G$ defined over $k$ which contains $t$.
Theorem~\ref{risp} implies that there exists an involution $n\in N(T)(\bar k)$ such that $nsn^{-1}=s^{-1}$ for all $s\in T$.
Thus $ntn=t^{-1}$ and $g\in n\Z_G(t)=nT$.
Let $g=ns_0$, for $s_0\in T$.
Then $g^2=ns_0ns_0=s_0^{-1}s_0=1$.
Hence $g$ is an involution and $g\in G(k)$.
Therefore $t$ is a product of two involutions $g$ and $gt$ in $G(k)$.

Suppose now $T$ is a maximal torus in $G$ defined over $k$ and $T(k)$ contains a strongly regular real element $t$.
Let $s\in T(k)$.
Suppose $g\in G(k)$ conjugates $t$ to $t^{-1}$.
Then we have proved that $g^2=1$.
We claim that $g$ conjugates $s$ to $s^{-1}$.
From calculations in the paragraph above, we have $g=ns_0$ for some $s_0\in T$.
Then $gsg^{-1}=ns_0ss_0^{-1}n^{-1}=nsn^{-1}=s^{-1}$.
But since $g$ is an involution in $G(k)$,  $s$ is a product of two involutions in $G(k)$. \qed

\vskip2mm
We note that in groups $G$ satisfying the hypothesis of the theorem, there are strongly regular elements in $G(k)$ which are not real in $G(k)$.
In~\cite{st} (see Theorem 6.3), it was shown that for a group $G$ of type $G_2$ defined over $k$, a semisimple element in $G(k)$ is real if and only if it is a product of two involutions in $G(k)$.
Examples of semisimple elements (\cite{st}, Theorem 6.10, Theorem 6.11 and Theorem 6.12) which are not real were also constructed in the same paper.
Hence in a maximal torus containing such an element no strongly regular element is real.

\subsection{An Obstruction to Reality}
The results in this subsection are known to experts (ref. \cite{s}, Section 11; \cite{se}, Chapter III, Section 2.3).
However, we include some with proofs for the sake of completeness.
Let $G$ be a connected linear algebraic group defined over a field $k$. 
In this section, we assume that the field $k$ is perfect.
We have, 
\begin{lemma}\label{realsemisimpleunipotent}
Let $g\in G$. 
Let $g=g_sg_u$ be the Jordan decomposition of $g$ in $G$. 
Let $H$ be the centralizer of $g_s$ in $G$.
Then, $g$ is real in $G$ if and only if $g_s$ is real and $g_u^{-1}, xg_ux^{-1}$ are conjugate in $H$, where $xg_sx^{-1}=g_s^{-1}$.
\end{lemma}
\noindent{\bf Proof :} Let $g$ be real in $G$, i.e., there exists $x\in G$ such that $xgx^{-1}=g^{-1}$.
Then $x$ conjugates $g_s$ and $g_u$ to $g_s^{-1}$ and $g_u^{-1}$ respectively.

Conversely let $h\in H$ such that $hg_u^{-1}h^{-1}=xg_ux^{-1}$.
Then, 
\begin{eqnarray*}
h^{-1}xg(h^{-1}x)^{-1}&=&h^{-1}xgx^{-1}h=h^{-1}xg_sx^{-1}xg_ux^{-1}h=h^{-1}g_s^{-1}xg_ux^{-1}h\\
&=&g_s^{-1}h^{-1}xg_ux^{-1}h=g_s^{-1}g_u^{-1}=g^{-1}.
\end{eqnarray*}
Hence $g$ is real in $G$.
 \qed

\vskip2mm
It is not true in general for an algebraic group $G$ that $g\in G$ is real if and only if $g_s$ is real and $g_u$ is real.
We give examples to illustrate this situation.

\vskip2mm
\noindent{\bf Example 1:} Let $G=GL_4(k)$.
We take $s=\diag(\lambda,\lambda,\lambda^{-1},\lambda^{-1})$ with $\lambda^2\neq 1$,\\\\ 
$u=\diag\left(\left(\begin{array}{cc}1&0\\ 0&1 \end{array}\right),\left(\begin{array}{cc}1&1\\ 0&1 \end{array}\right)\right)$
and $g=su$.
Then $g_s=s, g_u=u$ and the centralizer of $s$ in $G$ is $\Z_{GL_4(k)}(s)=\{\diag(A,B)\mid A,B\in GL_2(k)\}$.
The elements $s$ and $u$ are real but $g$ is not real.
In fact $xsx^{-1}=s^{-1}$ where 
$$x=\left(\begin{array}{cccc}0&0&1&0\\0&0& 0&1\\1&0&0&0\\0&1&0&0 \end{array}\right).$$ Any matrix 
$$y=\diag\left(\left(\begin{array}{cc}1&0\\ 0&1 \end{array}\right),\left(\begin{array}{cc}a&b\\ 0&-a \end{array}\right)\right)\in GL_4(k)$$ conjugates $u$ to $u^{-1}$.
The elements $u^{-1}$ and 
$$xux^{-1}=\diag\left(\left(\begin{array}{cc}1&1\\ 0&1 \end{array}\right),\left(\begin{array}{cc}1&0\\ 0&1 \end{array}\right)\right)$$ are not conjugate in $\Z_{GL_4(k)}(s)$. Hence $g$ is not real by Lemma~\ref{realsemisimpleunipotent}.

\vskip2mm
\noindent{\bf Example 2:} In $G=G_2$ over a finite field $k$, all semisimple as well as unipotent elements in $G(k)$ are strongly 
real but still there are non-real elements (ref. \cite{st}, Theorem 6.11).

\vskip2mm
Below we mention a cohomological obstruction to reality over the base field $k$.
Let $G$ be a connected linear algebraic group defined over $k$.
Suppose $t\in G(k)$ is real in $G(\bar k)$.
We put $H=\mathcal Z_G(t)$, the centralizer of $t$ in $G$.
Let $X=\{x\in G \mid xtx^{-1}=t^{-1}\}$.
Then $X$ is an $H$-torsor defined over $k$ with $H$-action given by $h.x=xh$ for $h\in H$ and $x\in X$.

Since $t$ is real over $\bar k$, we have $X\neq \phi$.
The torsor $X$ corresponds to an element of $H^1(k,H)$ (\cite{se}, Chapter 1, section 5.2, Proposition 33).
Let $x\in X$ and $\gamma$ be the cocycle corresponding to $X$.
Then $\gamma$ is given by $\gamma(\sigma)=x^{-1}\sigma(x)$ for all $\sigma\in \Gamma=\gal(\bar k/k)$.
We have,
\begin{proposition}
Let $G$ be a connected algebraic group defined over $k$.
Let $t\in G(k)$ be real over $\bar k$.
Then $t$ is real in $G(k)$ if and only $\gamma$, as above, represents the trivial cocycle in $H^1(k,H)$ where $H$ is the centralizer of $t$ in $G$.
\end{proposition}
\noindent{\bf Proof :}
Let $X$ be the $H$-torsor defined above. Then $\gamma\in H^1(k,H)$ is trivial if and only if $X$ has a $k$-rational point which is equivalent to $t$ is $k$-real.
\qed  

By the above, if $H^1(k,H)$ is trivial then $t$ is real in $G(k)$.
By a theorem of Steinberg (\cite{s}, Theorem 1.9; also see \cite{se}, Chapter III, section 2.3) if $H$ is a connected reductive group and $cd(k)\leq 1$ or $H$ is connected with $k$ perfect of $cd(k)\leq 1$, we have $H^1(k,H)=0$. 
In these situations $t$ is real.

\begin{proposition}\label{splittorusreal}
Let $G$ be a split connected semisimple adjoint group defined over an arbitrary field $k$ and suppose $-1$ belongs to the Weyl group of 
$G$. 
Let $T$ be a $k$-split maximal torus in $G$. 
Then every element of $T(k)$ is strongly real.
\end{proposition}
\noindent{\bf Proof :}
By Theorem~\ref{risp}, there exists $n_0\in N(T)(\bar k)$ such that ${n_0}^2=1$ and $n_0s{n_0}^{-1}=s^{-1}$ for all $s\in T$. Consider the Galois cocycle $\gamma(\sigma)=n_0\sigma(n_0)$ for $\sigma\in\Gamma=\gal(\bar k/k)$. Since $T$ is defined over $k$, we have, for $s\in T$ and $\sigma\in\Gamma$,
$$\sigma(n_0)s\sigma(n_0)^{-1}=\sigma(n_0\sigma^{-1}(s)n_0)=\sigma(\sigma^{-1}(s^{-1}))=s^{-1}.$$
Hence, we must have, in the Weyl group $W=N(T)/T$, $n_0T=\sigma(n_0)T$. Therefore $\gamma(\sigma)=n_0\sigma(n_0)\in T$. Hence $\gamma$ is 
a $1$-cocycle in $H^1(k,T)$. But since $T$ is $k$-split, $H^1(k,T)=0$. Hence there is $s\in T$ such that 
$$\gamma(\sigma)=n_0\sigma(n_0)=s^{-1}\sigma(s).$$
This gives $sn_0=\sigma(sn_0)$ for all $\sigma\in\Gamma$. Hence $sn_0\in T(k)$. Also 
$$(sn_0)^2=sn_0sn_0=ss^{-1}=1.$$ 
Therefore $g=sn_0$ is an involution in $T(k)$ and for any $t\in T(k)$, we have,
$$ gtg^{-1}=gtg=sn_0tn_0s^{-1}=st^{-1}s^{-1}=t^{-1}.$$
Thus $(gt)^2=1$ and $t=g.gt$. Hence $t$ is strongly real.
\qed

\subsection{Reality over Fields of $cd(k)\leq 1$}
In this section we discuss reality for algebraic groups over fields $k$ with $cd(k)\leq 1$.
We have,
\begin{theorem}\label{semreal}
Let $k$ be a field with $cd(k)\leq 1$. Let $G$ be a connected reductive group defined over $k$ with $-1$ in its Weyl group. Then every semisimple element in $G(k)$ is real in $G(k)$.
\end{theorem}
\noindent{\bf Proof :}
Let $t\in G(k)$ be semisimple. Let $T$ be a maximal torus defined over $k$ with $t\in T(k)$. Let $W=N(T)/T$ be the Weyl group of $G$, 
where $N(T)$ is the normalizer of $T$ in $G$. We have the exact sequence
$$1\rightarrow T\rightarrow N(T)\rightarrow W\rightarrow 1.$$The corresponding Galois cohomology sequence is
$$1\rightarrow T(k)\rightarrow N(T)(k)\rightarrow W(k)\rightarrow H^1(k,T)\rightarrow\cdots .$$
Since $cd(k)\leq 1$, by Steinberg's theorem (\cite{s}, Theorem 1.9), $H^1(k,T)=0$. Hence the longest element $w_0$ in the Weyl group, which acts by $-1$ on the set of roots, lifts to an element $h\in N(T)(k)$. Hence $hth^{-1}=t^{-1}$ with $h\in G(k)$ and $t$ is real in $G(k)$.
\qed  
\begin{corollary}
Let $G$ and $k$ be as in the above theorem. Then every regular element of $G$ is real.
\end{corollary}
\noindent{\bf Proof :}
Let $g\in G$ be regular and $g=g_sg_u$ be the Jordan decomposition of $g$ in $G$ with $g_s$ semisimple and $g_u$ unipotent. Then, by the above theorem, $hg_sh^{-1}=g_s^{-1}$ for some $h\in G$. Then $hg_uh^{-1}$ and $g_u^{-1}$ are regular unipotents in $\mathcal Z_G(g_s)^0$ and hence there is $x\in\mathcal Z_G(g_s)$ such that $xhg_uh^{-1}x^{-1}=g_u^{-1}$. Then $(xh)g(xh)^{-1}=g^{-1}$ and hence $g$ is real (see Corollary 1.9, Chapter III, \cite{ss}). \qed

\begin{theorem}\label{regreal}
Let $k$ be a field with $cd(k)\leq 1$. 
Let $G$ be a connected semisimple adjoint group defined over $k$ with $-1$ in its Weyl group.
Then every semisimple element in $G(k)$ is strongly real in  $G(k)$.
\end{theorem}
\noindent{\bf Proof : }
We may assume $G$ is simple. 
Let $t\in G(k)$ be a semisimple element.
Let $T$ be a maximal torus in $G$ defined over $k$ which contains $t$, i.e., $t\in T(k)$. Since $cd(k)\leq 1$, by Steinberg (\cite{s}, Theorem 1.9) we have $H^{1}(k,T)=0$. The rest of the proof follows exactly along the lines of the proof of Proposition~\ref{splittorusreal}. 
\qed

\noindent
{\bf Remark :} It seems likely that the results of this section are valid over non-perfect fields also, however we have not been able to prove this.  
\section{Reality in Classical Groups}
In this section we discuss structure of real elements in  classical groups.
We assume $k$ is an arbitrary field of characteristic not $2$.

\subsection{The Groups $GL_n(k)$ and $SL_n(k)$}

It was proved by Wonenburger (\cite{w1}, Theorem 1) that an element of $GL_n(k)$ is real if and only if it is strongly real in $GL_n(k)$.
However, similar result is false for matrices over division algebras.
In~\cite{e}, (Lemma 2 and Lemma 3) Ellers constructs an example  of a simple transformation on a vector space $V$ over the real quaternion division algebra $\mathbb H$, which is conjugate to its inverse but is not a product of two involutions.
This is also evident by looking at the following example.
Let $\mathbb H=\mathbb R.1\oplus \mathbb R.i\oplus \mathbb R.j\oplus \mathbb R.ij$ where $i,j,k$ have usual meanings.
In the group $GL_1(\mathbb H)$, the element $i$ is conjugate to its inverse by $j$.
The only nontrivial element of $GL_1(\mathbb H)$ which is an involution is $-1$ and hence $i$ is not a product of two involutions in $GL_1(\mathbb H)$.

In this section, we explore the structure of real elements in $SL_n(k)$.
We follow the proof of Wonenburger for $GL_n(k)$ (ref. \cite{w1}, Theorem 1) and modify it for our purpose.
\begin{theorem}\label{realsln}
Let $V$ be a vector space of dimension $n$ over $k$ and let $t\in SL(V)(k)$.
Suppose $n\not\equiv 2 \imod 4$.
Then $t$ is real in $SL(V)(k)$ if and only if $t$ is strongly real in $SL(V)(k)$.
\end{theorem}
\noindent{\bf Proof :} 
Let $\delta_1(X),\ldots,\delta_n(X)$ be the invariant factors of $t$ in $k[X]$.
Since $t$ is real, each $\delta_i(X)$ is self-reciprocal.
The space $V$ decomposes as $V=\oplus_{i=1}^n V_i$, where each $V_i$ is a cyclic, $t$ invariant subspace of $V$ and the minimal polynomial of $t_i=t|_{V_i}$ is the self-reciprocal polynomial $\delta_i(X)$.
We shall construct involutions $H_i$ in $GL(V_i)$, conjugating $t_i$ to $t_i^{-1}$, with $\det(H_i)=(-1)^m$ if dimension of $V_i=2m$ and $\det(H_i)=(-1)^m$ or $(-1)^{m+1}$ when dimension of $V_i=2m+1$.
Then $H=\oplus_{i=1}^nH_i$ is an involution conjugating $t$ to $t^{-1}$ and $\det(H)=1$ if $\dim(V)\not\equiv 2 \imod 4$.

Now $t_i$ is a cyclic linear transformation on the vector space $V_i$ with self-reciprocal characteristic polynomial $\chi_{t_i}(X)=\delta_i(X)$. Hence, we can write $\chi_{t_i}(X)=(X-1)^r(X+1)^sf(X)$ where $f(\pm 1)\neq 0$ and $V_i=W_{-1}\oplus W_{1}\oplus W_0$, where $W_{-1}, W_1$ and $W_0$ are the kernels of $(t_i-1)^r, (t_i+1)^s$ and $f(t_i)$ respectively.
To produce the involution $H_i$ on $V_i$ as above, it suffices to do so on each of $W_{-1}, W_1$ and $W_0$.
Hence we are reduced to the following cases.
Let $t$ be a cyclic linear transformation on a vector space $V$ with self reciprocal characteristic polynomial $\chi_{t}(X)$, of the following two types;
\begin{enumerate}
\item the degree of $\chi_{t}(X)$ is even, say $2m$,
\item $\chi_{t}(X)=(X-1)^{2m+1}$ or $(X+1)^{2m+1}$.
\end{enumerate}
We claim that in the first case $t$ is conjugate to $t^{-1}$ by an involution whose determinant is $(-1)^m$ and in the second, 
there are involutions with determinant $(-1)^{m}$ or $(-1)^{m+1}$ conjugating $t$ to $t^{-1}$.
We first prove that in both the cases, $V$ admits a decomposition $V=V_+\oplus V_-$, 
invariant under $t+t^{-1}$ and such that $(t-t^{-1})V_{\pm}\subset V_{\mp}$.

In the first case, since $V$ is cyclic, there is a vector $u\in V$ such that $\mathcal E=\{u,tu,\ldots,t^{2m-1}u\}$ is a basis of $V$. 
Set $S^mu=y$. Then   
$$
\mathcal B=\{y,(t+t^{-1})y,\ldots,(t^{m-1}+t^{-m+1})y,(t-t^{-1})y,\ldots,(t^m-t^{-m})y\}
$$
is a basis of $V$. 
Let $V_+$ denote the subspace generated by the first $m$ vectors of $\mathcal B$ and $V_-$ that by the latter $m$ vectors.
Then $t+t^{-1}$ leaves $V_+$ as well as $V_-$ invariant, $(t-t^{-1})V_{\pm}\subset V_{\mp}$ and $V=V_+\oplus V_-$. 
In the second case, we take 
$$  
\mathcal B=\{y,(t+t^{-1})y,\ldots,(t^{m}+t^{-m})y,(t-t^{-1})y,\ldots,(t^m-t^{-m})y\}
$$
as a basis of $V$ and $V_+$ as the span of the first $m+1$ vectors from $\mathcal B$ and $V_-$ as the span of the latter $m$.
In the first case, let $H=1|_{P}\oplus -1|_{Q}$.
Then $H$ is an involution which conjugates $t$ to $t^{-1}$ and has determinant $(-1)^m$.
In the second case, we consider $H_1=1|_{V_+}\oplus -1|_{V_-}$ and $H_2=-1|_{V_+}\oplus 1|_{V_-}$.
Then $H_1$ and $H_2$ both are involutions which conjugate $t$ to $t^{-1}$ and have determinants $(-1)^m$ and $(-1)^{m+1}$ respectively.   \qed

\vskip2mm
\noindent{\bf Remarks :} 

{\bf 1.} An element $S=\diag(\alpha,\alpha^{-1},\beta,\beta^{-1},\gamma,\gamma^{-1})\in SL_6(k)$ such that all the diagonal entries are distinct, can be conjugated to its inverse by  
$$
H=\diag\left(\left (\begin{array}{cc} 0&-1 \\  1&0 \\ \end{array}  \right),\left (\begin{array}{cc} 0&-1 \\  1&0 \\ \end{array}  \right),\left (\begin{array}{cc} 0&-1 \\  1&0 \\ \end{array}  \right) \right)\in SL_6(k)
$$ 
where $H^2=-1$. 
In fact any element $T\in SL_6(k)$ such that $TST^{-1}=S^{-1}$ is of the form:
$$
T=\diag\left(\left (\begin{array}{cc} 0&a \\  \tilde a&0 \\ \end{array}  \right),\left (\begin{array}{cc} 0&b \\  \tilde b&0 \\ \end{array}  \right),\left (\begin{array}{cc} 0&c \\ \tilde c&0 \\ \end{array}  \right) \right)
$$ 
where $a\tilde ab\tilde bc\tilde c=-1$.
Suppose $T^2=1$.
Then $a\tilde a=1, b\tilde b=1, c\tilde c=1$.
This implies that $a\tilde ab\tilde bc\tilde c =1$, a contradiction.
Hence there is no involution in $SL_6(k)$ conjugating $S$ to $S^{-1}$, i.e., $S$ is real semisimple but not strongly real in $SL_6(k)$.\\

{\bf 2.} Let us take $A=\left (\begin{array}{cc} 1&1 \\  0&1 \\ \end{array}  \right) $, a unipotent element in $SL_2(k)$.
Then any element $X\in GL_2(k)$ such that $XAX^{-1}=A^{-1}$ has the form $X=\left (\begin{array}{cc} a&b \\  0&-a \\ \end{array}  \right) $.
Then, $A$ is conjugate to $A^{-1}$ in $SL_2(k)$ if and only if $-1$ is a square in $k$.
In that case ($-1$ is a square in $k$) the element $X$ which conjugates $A$ to its inverse satisfies $X^2=-1$, not an involution, and hence $A$ is not strongly real in $SL_2(k)$.

\subsection{Groups of Type $A_1$}
In this section we study real semisimple elements in $SL_2(k)$ and $PSL_2(k)=SL_2(k)/\Z(SL_2(k))$.
Though the proofs of Corollary~\ref{realpsl2kbar}, Proposition~\ref{realpsp2} and \ref{gptypesl2} follow essentially from 
Theorem~\ref{streal}, we give proofs with explicit computations.
We fix an algebraic closure $\bar k$ of $k$.
Let $G=SL_2(\bar k)$. 
Fix the maximal torus $T=\{\diag(\alpha,\alpha^{-1})\mid \alpha\in\bar k^*\}$ in $G$.
\begin{lemma}\label{realsl2kbar}
Every semisimple element of $G=SL_2(\bar k)$ is real in $G$.
The only involutions in $G$ are $\{1,-1\}$, hence non-central semisimple elements are not a product of involutions in $G$.
Moreover, every semisimple element of $G$ is conjugate to its inverse by an involution in $GL_2(\bar k)$, hence is strongly real in $GL_2(\bar k)$.
\end{lemma}
\noindent{\bf Proof : } Let $t\in SL_2(\bar k)$ be semisimple.\\
First, assume that $t=\diag(\alpha,\alpha^{-1})\in T$.
Let $g=\left (\begin{array}{cc} 0&-1 \\  1&0 \\ \end{array}  \right) \in SL_2(\bar k)$.
Then $g^2=-1$ and $gtg^{-1}=t^{-1}$.
Hence, for any $t\in T$, $gtg^{-1}=t^{-1}$.

Now let $n=\left (\begin{array}{cc} 0&1 \\  1&0 \\ \end{array}  \right)$.
Then we have, for any $t\in T$, $ntn^{-1}=t^{-1}$ and $n$ is an involution with $\det(n)=-1$. 
Hence, for any $t\in T$, we have $t=n.nt$, a product of two involutions in $GL_2(\bar k)$.
If $s\in SL_2(\bar k)$ is semisimple then $gsg^{-1}\in T$ for some $g\in SL_2(\bar k)$.
If $gsg^{-1}=\rho_1\rho_2$, $\rho_i\in GL_2(\bar k), \rho_i^2=1$, then $s=g^{-1}\rho_1g.g^{-1}\rho_2g$, and $g^{-1}\rho_ig$ are involutions in $GL_2(\bar k)$.  \qed

\begin{corollary}\label{realpsl2kbar}
Let $G=PSL_2(\bar k)$ and $t$ be a semisimple element in $G$.
Then $t$ is real in $G$ if and only if $t$ is strongly real in $G$.
\end{corollary}
\noindent{\bf Proof :} Let $t$ as above be real. 
Let $t_0\in SL_2(\bar k)$ be a representative of $t$.
Then $t_0$ is either conjugate to $t_0^{-1}$ or $-t_0^{-1}$ in $SL_2(\bar k)$.
When $t_0$ is conjugate to $t_0^{-1}$, it follows from the previous lemma that there exists an element $s\in SL_2(\bar k)$ with $s^2=-1$ such that $st_0s^{-1}=t_0^{-1}$.
We have $t_0=(-s).(st_0)$ and hence $t$ as a product of two involutions in $PSL_2(\bar k)$.

Now suppose $t_0$ is conjugate to $-t_0^{-1}$ in $SL_2(\bar k)$.
Then the characteristic polynomial of $t_0$ is $X^2+1$.
In this case $t$ itself is an involution in $PSL_2(\bar k)$.    \qed

\vskip 2mm
We need,.
\begin{lemma}\label{stregularcentral}
Let $t\in SL_2(k)$ be a semisimple element.
Then $t$ is either strongly regular or central in $SL_2(k)$.
\end{lemma}

\vskip2mm
Hence we can produce real elements in $SL_2(k)$, as in Lemma \ref{realsl2kbar}, which are not a product of two involutions in $SL_2(k)$.
\begin{proposition}\label{realpsp2}
Let $t\in PSL_2(k)$ be a semisimple element.
Then $t$ is real in $PSL_2(k)$ if and only if $t$ is strongly real in $PSL_2(k)$.
\end{proposition}
\noindent{\bf Proof : } Let $t_0\in SL_2(k)$ be a representative of $t$.
Since $t$ is real in $PSL_2(k)$, it follows that $t_0$ is either conjugate to ${t_0}^{-1}$ or ${-t_0}^{-1}$ in $SL_2(k)$.
In the second case, the characteristic polynomial of $t_0$ must be $X^2+1$ and hence ${t_0}^2=-1$.
For the first case we prove that there exists $s \in SL_2(k)$ with $s^2=-1$ such that $st_0s^{-1}=t^{-1}$.

If $t_0$ is central, it is either $1$ or $-1$.
Hence we may assume that the element $t_0$ is conjugate to the matrix $t_1=\diag(\alpha,\alpha^{-1})$ in $SL_2(\bar k)$, for some $\alpha\in \bar k$ with $\alpha^2\neq 1$.
Let $$n=\left(\begin{array}{cc} 
0&-1\\1&0
\end{array}\right)\in SL_2(\bar k).$$
Then $nt_1n^{-1}=t_1^{-1}$ and $n^2=-1$.
In fact $n$ conjugates every element of the torus $T_1=\{\diag(\gamma,\gamma^{-1})\mid \gamma\in \bar k^*\}$ to its inverse.
Hence there exists $h\in SL_2(\bar k)$ such that $ht_0h^{-1}={t_0}^{-1}$ and $h^2=-1$.
Moreover, $h$ conjugates every element of the maximal torus $T$ containing $t_0$, to its inverse. 
Since $t_0$ is real in $SL_2(k)$, there exists $g\in SL_2(k)$ such that $gt_0g^{-1}={t_0}^{-1}$.
Then $g\in h\Z_{SL_2(\bar k)}(t_0)$.
Since $t_0$ is not central (by Lemma~\ref{stregularcentral}) we have $\Z_{SL_2(\bar k)}(t_0)=T$.
We write $g=hx$ where $x\in T$.
Then $g^2=hxhx=-hxh^{-1}x=-x^{-1}x=-1$ and this proves the required result. 
\qed

\vskip2mm
We now consider $Q$, a quaternion algebra over $k$.
It is a central simple algebra over $k$ of degree $2$.
We note that $SL_1(Q)=\{x\in Q^*\mid Nrd(x)=1\}$ is a form of $SL_2$ over $k$.
We denote the group $SL_1(Q)/\Z(SL_1(Q))$ by $PSL_1(Q)$.
\begin{proposition}\label{gptypesl2}
Let $G=PSL_1(Q)$ and $t\in G$ be a semisimple element.
Then, $t$ is real in $PSL_1(Q)$ if and only if $t$ is strongly real in $PSL_1(Q)$.
Furthermore, $G=SL_1(Q)$ has real elements which are not strongly real.
\end{proposition}
\noindent{\bf Proof :} We first observe that an element $t\in Q^*$ is either strongly regular or central. 
Proof of this fact and the rest of the proposition is on similar lines as in Lemma~\ref{stregularcentral} and Proposition~\ref{realpsp2}. \qed

\subsection{$SL_1(D)$, deg(D) Odd}
We now consider anisotropic simple groups of type $A_n$, for $n$ even.
These are the groups $SL_1(D)$ for central division algebras of degree $n+1$.
Let $D$ be a central division algebra over a field $k$, with degree $D$ odd.
Let $G=D^*$ or $G=SL_1(D)=\{x\in D^*\mid Nrd(x)=1\}$.
We have,
\begin{theorem}
Let $G$ be as above.
Then the only real elements in $G=D^*$ are $\pm 1$.
In $G=SL_1(D)$, there are no nontrivial real elements.
\end{theorem}
\noindent{\bf Proof :}
We first prove that there are no non-central real element in $G$ and there are no non-central involutions in $G$.
Let $t\in G$ be a real element which is not in the center of $D$.
Then $k(t)$ is a subfield $\neq k$ contained in $D$ and has a field automorphism given by $t\mapsto t^{-1}$ of order two.
Hence the degree of $k(t)$ over $k$ is even. 
But degree of $D$ being odd, $D$ can not contain a field extension of even degree. 
Hence there are no real elements which are not in the center of $G$. 

Now let $t\in G$ be a non-central involution.
Then $k(t)$ is a field extension over $k$ of even degree.
Following similar argument as in the previous paragraph, we get a contradiction.
Hence any involution in $G$ is in the center of $G$. 
Since $D$ is central and degree $D$ is odd, any such involution is trivial.
This completes the proof. \qed

\begin{corollary}\label{formorthogonalgroup}
Let $D$ be a central division algebra over a field $k$, with degree $D$ odd.
Let $\sigma$ be an involution on $D$.
Then the group $Iso(D,\sigma)=\{x\in D\mid x\sigma(x)=1\}$ has no nontrivial real elements.
\end{corollary}
\noindent {\bf Proof :}
Since $Iso(D,\sigma)\subset D^*$, the result follows from the above theorem.
\qed

We remark that (\cite{kmrt}, Corollary 2.8 and Section 12.B) the group $Iso(D,\sigma)$, for $\sigma$ of the first kind, is a form of the orthogonal group.
The group $Iso(D,\sigma)$, for $\sigma$ of the second kind, is a form of the unitary group. Hence the results above prove the absence of nontrivial real elements in anisotropic $k$-forms of orthogonal and unitary groups when the degree of the underlying division algebra is odd.

\subsection{Orthogonal Groups}
Let $V$ be a vector space over $k$ with a nondegenerate quadratic form $Q$.
We denote the orthogonal group by $O(Q)$.
Then Wonenburger proved (\cite{w1}, Theorem 2),
\begin{proposition}
 Any element of the orthogonal group $O(Q)$ is strongly real, i.e., the group $O(Q)$ is bireflectional.
Hence every element of $O(Q)$ is real.
\end{proposition}
\noindent Djokovic (\cite{d}, Theorem 1) extended this result to fields of characteristic $2$.
However Knuppel and Nielsen proved (\cite{kn}, Theorem A),
\begin{proposition}\label{theoremkn1}
The group $SO(Q)$ is trireflectional, except when $\dim(V)=2$ and $V\not= \mathcal H_3$, where $\mathcal H_3$ is the hyperbolic plane over $\mathbb F_3$.
The group $SO(Q)$ is bireflectional if and only if $\dim(V)\not\equiv 2 \imod 4$ or $V=\mathcal H_3$, and hence in this case every element is real.
\end{proposition}
\noindent They give necessary and sufficient condition for an element in special orthogonal group to be strongly real (\cite{kn}, Proposition 3.3).
\begin{proposition}\label{theoremkn2}
Let $t\in SO(Q)$.
Then $t$ is a product of two involutions in $SO(Q)$ if and only if $\dim(V)\not\equiv 2\imod 4$ or an orthogonal decomposition of $V$ into orthogonally indecomposable $t$-modules contains an odd dimensional summand. 
\end{proposition}
\noindent {\bf Proof :}
We shall indicate the proof when $t$ is semisimple, since that concerns us. 
Note that when $\dim(V)=2$, any $\rho\in O(Q)-SO(Q)$ satisfies $\rho^2=1$ and 
$\rho t\rho^{-1}=t^{-1}$. Let $t\in SO(Q)$ be any semisimple element, where $\dim(V)\not\equiv 2\imod 4$. 
Let $\bar V=V\otimes \bar k$ and, for $\alpha\in {\bar k}^*$, let $\bar V_{\alpha}=\{x\in\bar V|t(x)=\alpha x\}$ and 
$\bar W_{\alpha}=\bar V_{\alpha}\oplus\bar V_{\alpha^{-1}}$. Then $\bar W_{\alpha}$ is nondegenerate and defined over the subfield $k_{\alpha}$ of $\bar k$ which is the fixed field of the subgroup of $\Gamma=\gal(\bar k/k)$ fixing the unordered pair $\{\alpha,\alpha^{-1}\}$. Let $W_{\alpha}$ denote the descent of $\bar W_{\alpha}$ over $k_{\alpha}$. Then $W_{\alpha}$ is a direct sum of $m_{\alpha}$ (say) 
$2$-dimensional subspaces, on each of which ( a conjugate of ) $t$ restricts to  $\diag\{\alpha,\alpha^{-1}\}\in SO(W_{\alpha})$. Then by the $2$-dimensional situation, there is $g_{\alpha}\in O(W_{\alpha})-SO(W_{\alpha})$, such that $g_{\alpha}^2=1$ and 
$g_{\alpha}tg_{\alpha}^{-1}=t^{-1}$. Let $W_{\Gamma\alpha}=\oplus_{\sigma\in\Gamma}W_{\sigma\alpha}$ and $g_{\Gamma\alpha}=\oplus_{\sigma\in\Gamma}g_{\sigma\alpha}$. Then $W_{\Gamma\alpha}$ and $g_{\Gamma\alpha}$ are defined over $k$, $g_{\Gamma\alpha}^2=1$ and 
$g_{\Gamma\alpha}tg_{\Gamma\alpha}^{-1}=t^{-1}$ on $W_{\Gamma\alpha}$. Since $V$ is the orthogonal direct sum of $V_{\pm 1}$ and the subspaces $W_{\Gamma\alpha}$, the result follows from the fact that the determinant of $g_{\Gamma\alpha}=(-1)^{\frac{1}{2}\dim W_{\Gamma\alpha}}$.
\qed        

Now we take up the case $\dim(V)\equiv 2\imod 4$.
First we prove,
\begin{lemma}\label{so6real}
Let $t\in SO(Q)$ where $\dim(V)\equiv 2\imod 4$.
Let $t$ be a semisimple element which has only two distinct eigenvalues $\lambda$ and $\lambda^{-1} ($hence $\lambda\not=\pm 1)$ over $\bar k$.
Then $t$ is not real in $SO(Q)$.
\end{lemma}
\noindent{\bf Proof :} 
We prove that the element $t$ is not real over $\bar k$.
Let $\dim(V)=2m$ where $m$ is odd.
The element $t$ over $\bar k$ is conjugate to $A=\diag(\underbrace{\lambda,\ldots,\lambda}_m,\underbrace{\lambda^{-1},\ldots,\lambda^{-1}}_m)$ with $\lambda\neq\pm 1$ in $SO(J)$ where $J$ is the matrix of the quadratic form over $\bar k$ given by $J=\left(\begin{array}{cc} 0&S\\S&0 \end{array}\right)$ where 
$$S=\left(\begin{array}{ccccc} 0&0&\ldots&0&1\\0&0&\ldots&1&0\\\vdots&&&&\vdots\\1&0&\ldots&0&0 \end{array}\right),$$ an $m\cross m$ matrix.
Now suppose $A$ is real in $SO(J)$, i.e., there exists $T\in SO(J)$ such that $TAT^{-1}=A^{-1}$.
Then $T$ maps the $\lambda$-eigen subspace of $A$ to the $\lambda^{-1}$-eigen subspace of $A$ and vice-versa.
Hence $T$ has the following form:
$$T=\left(\begin{array}{cc} 
0&B\\C&0
\end{array}\right)$$
for $m\cross m$ matrices $B$ and $C$.
Since $T$ is orthogonal, it satisfies $\tra TJT=J$, which gives $\tr BSC=S$.
That is, $\det(B)\det(C)=1$.
Hence $\det(T)=(-1)^{m}\det(B)\det(C)=-\det(B)\det(C)=-1$ since $m$ is odd.
This contradicts that $T\in SO(J)$.
Hence $A$ is not real in $SO(J)$ and hence $t$ is not real in $SO(Q)$.             \qed

\begin{lemma}
Let $\dim(V)\equiv 0\imod 4$ and $t\in SO(Q)$ be semisimple.
Suppose $t$ has only two distinct eigenvalues $\lambda$ and $\lambda^{-1} ($hence $\lambda\not=\pm 1)$ over $\bar k$.
Then, any element $g\in O(Q)$ such that $gtg^{-1}=t^{-1}$ belongs to $SO(Q)$, i.e., $\det(g)=1$. 
\end{lemma}
\noindent{\bf Proof :}
We follow the notation in the previous lemma.
Let $\dim(V)=2m$, where $m$ is even.
As in the proof of the previous lemma, we may assume $t$ is diagonal.
Then any element $T$ that conjugates $t$ to $t^{-1}$ over $\bar k$, is of the form $T=\left(\begin{array}{cc} 
0&B\\C&0
\end{array}\right)$. 
We have $\det(T)=(-1)^{m}\det(B)\det(C)=\det(B)\det(C)=1$. 
Since $g$ is a conjugate of $T$, the claim follows.       \qed
\vskip2mm

\noindent Now we state the main theorem about special orthogonal groups.
\begin{theorem}\label{realsoq}
Let $Q$ be a nondegenerate quadratic form on $V$, with dimension of $V$ arbitrary. Let $t\in SO(Q)$ be a semisimple element.
Then, $t$ is real in $SO(Q)$ if and only if $t$ is strongly real in $SO(Q)$.
\end{theorem}
\noindent{\bf Proof :}
If $\dim(V)\not\equiv 2\imod 4$ then the theorem follows from Propositions~\ref{theoremkn1} and \ref{theoremkn2}. 
Hence let us assume that $\dim(V)\equiv 2\imod 4$.
Let $\dim(V)=2m$ where $m$ is odd.
In this case we will prove that the element $t$ is real in $SO(Q)$ if and only if $1$ or $-1$ is an eigenvalue of $t$.

First we prove that if $1$ and $-1$ are not eigenvalues then $t$ is not real.
It is enough to prove this statement over $\bar k$.
We write $\bar V=V\tensor_k \bar k$ and continue to denote $t$ over $\bar k$ by $t$ itself.
We have a $t$-invariant orthogonal decomposition of $\bar V$;  
$$\bar V=\bar V_1\oplus \bar V_{-1}\oplus \bar V_{\lambda_1^{\pm 1}}\oplus \ldots\oplus \bar V_{\lambda_r^{\pm 1}}$$ 
where $\bar V_1$ and $\bar V_{-1}$ are the eigenspaces of $t$ corresponding to $1$ and $-1$ respectively and $\bar V_{\lambda_j^{\pm 1}}=\bar V_{\lambda_j}\oplus \bar V_{\lambda_j^{-1}}$  where $\bar V_{\lambda_j}$ is the eigenspace corresponding to $\lambda_j$ for $\lambda_j^2\neq 1$.
Since $1$ and $-1$ are not eigenvalues for $t$, we have $\bar V_1=0$ and  $\bar V_{-1}=0$.
If $r=1$ it follows from Lemma~\ref{so6real} that $t$ is not real.
Hence we may assume $r\geq 2$.
We denote the restriction of $t$ on $\bar V_{\lambda_j^{\pm 1}}$ by $t_j$.
Let the dimension of $\bar V_{\lambda_j^{\pm 1}}$ be $n_j$.
Since $\lambda_j \neq \pm 1$, $n_j$ is even and is either $0\imod 4$ or $2\imod 4$.
Let the number of subspaces $\bar V_{\lambda_j^{\pm 1}}$ such that $n_j$ is $2\imod 4$ be $s$.
Then $s$ is odd, since $\dim(V)\equiv 2\imod 4$.
Let $g\in SO(Q)$ such that $gtg^{-1}=t^{-1}$.
Then $g$ leaves $\bar V_{\lambda_j^{\pm 1}}$ invariant for all $j$.
We denote the restriction of $g$ on $\bar V_{\lambda_j^{\pm 1}}$ by $g_j$.
Then $g_j\in O(\bar V_{\lambda_j^{\pm 1}})$ and $g_jt_jg_j^{-1}=t_j^{-1}$.
From the previous lemma, determinant of $g_j$ is $1$ whenever $n_j\equiv 0\imod 4$ and the determinant of $g_j$ is $-1$ whenever $n_j\equiv 2\imod 4$.
Hence the determinant of $g$ is $(-1)^{s}=-1$, which contradicts $g\in SO(Q)$.
Hence $t$ can not be real in $SO(Q)$.

Conversely, if $1$ or $-1$ is an eigenvalue then the subspace $\bar V_1$ or $\bar V_{-1}$ is non-zero.
These subspaces are defined over $k$.
Let us denote their descents by $V_1$ and $V_{-1}$ over $k$.
Then the dimension of $V_1$ and $V_{-1}$ is even, since $\dim(V)\equiv 2\imod 4$. Note that the restrictions of $t$ to $V_1$ and $V_{-1}$ 
are respectively $1$ and $-1$. Write the restriction of $t$ to $W=(V_1\oplus V_{-1})^{\perp}$ as a product of two involutions in $O(W)$. If any of these involutions has determinant $-1$, we write $1$ and $-1$ respectively on $V_1$ and $V_{-1}$ as a product of two involutions, 
each having determinant $1$ or $-1$, 
adjusted suitably, so as to get an expression of $t$ as a product of two involutions in $SO(Q)$. 
\qed

\subsection{Symplectic Groups}

Now we consider the symplectic group.
Let $V$ be a vector space of dimension $2n$ with a nondegenerate symplectic form.
We denote the symplectic group by $Sp(2n,k)$.
The center of this group is $\Z(Sp(2n,k))=\{\pm 1\}$. We denote the projective symplectic group by $PSp(2n,k)=Sp(2n,k)/\Z(Sp(2n,k))$. We begin by proving results for reality in $PSp(2,\bar k)$ and $PSp(4,\bar k)$, which we use for the general case.
\begin{lemma}\label{sp2kbar}
Let $t\in Sp(2,\bar k)$ be a semisimple element.
Suppose that $t$ is either conjugate to $t^{-1}$ or $-t^{-1}$.
Then the conjugation can be achieved by an element $s\in Sp(2,\bar k)$ such that $s^2=-1$.
Hence a semisimple element of $PSp(2,\bar k)$ is real if and only if it is strongly real in $PSp(2,\bar k)$.
\end{lemma}
\noindent{\bf Proof :} We note that $Sp(2,\bar k)=SL(2,\bar k)$. 
Hence proof follows from Corollary~\ref{realpsl2kbar}.
\qed

\begin{lemma}\label{sp4kbar}
Let $t\in Sp(4,\bar k)$ be a semisimple element.
Suppose that $t$ is either conjugate to $t^{-1}$ or $-t^{-1}$.
Then the conjugation can be achieved by an element $s\in Sp(4,\bar k)$ such that $s^2=-1$.
Hence a semisimple element of $PSp(4,\bar k)$ is real if and only if it is strongly real in $PSp(4,\bar k)$.
\end{lemma}
\noindent{\bf Proof :}  Let $J=\diag\left(\left (\begin{array}{cc} 0&-1 \\  1&0 \\ \end{array}  \right),\left (\begin{array}{cc} 0&-1 \\  1&0 \\ \end{array}  \right) \right)$.
Then $Sp(4,\bar k)=\{A\in GL(4,\bar k)\mid \tr AJA=J\}$.
We first assume $t$ is conjugate to $t^{-1}$.
We may assume $t=\diag(\lambda,\lambda^{-1},\mu,\mu^{-1})$.
We let 
$$
g=\diag\left(\left (\begin{array}{cc} 0&-1 \\ 1&0 \\ \end{array}  \right),\left (\begin{array}{cc} 0&-1 \\  1&0 \\ \end{array}  \right) \right) \in Sp(4,\bar k).
$$ 
Then $g^2=-1$ and $gtg^{-1}=t^{-1}$.

Now let $t$ be conjugate to $-t^{-1}$.
Then we may assume $t=\diag(\lambda,\lambda^{-1},-\lambda,-\lambda^{-1})$.
Let $$g=\left (\begin{array}{cccc} 0&0&0&-1 \\ 0&0&1&0 \\0&-1&0&0\\1&0&0&0 \end{array}  \right).$$
Then $g$ belongs to $Sp(4,\bar k)$ with $g^2=-1$ and $gtg^{-1}=-t^{-1}$.
\qed

\begin{theorem}\label{realpsp}
Let $t\in Sp(2n,k)$ be a  semisimple element.
Suppose $t$ is either conjugate to $t^{-1}$ or $-t^{-1}$.
Then the conjugation can be achieved by an element $s\in Sp(2n,k)$ such that $s^2=-1$.
Hence a semisimple element of $PSp(2n,k)$ is real if and only if it is strongly real in $PSp(2n,k)$.
\end{theorem}
\noindent{\bf Proof : } 
First we consider semisimple elements in $Sp(2n,\bar k)$.
Let $t\in Sp(2n,\bar k)$ be semisimple with $t$ conjugate to $t^{-1}$.
Then $t$ can be conjugated to 
$$\diag(\lambda_1,\lambda_1^{-1},\ldots,\lambda_n,\lambda_n^{-1})$$
 and this diagonal element can be conjugated to its inverse by $s=\diag(\underbrace{N,\ldots,N}_n )$
where  $N=\left(\begin{array}{cc}
0&-1\\ 1&0
\end{array}\right)$. 
Clearly $s^2=-1$. 
A conjugate of $s$ then does the job.

Now let us assume $t$ is conjugate to $-t^{-1}$ in $Sp(2n,\bar k)$.
Then $t$ can be conjugated to $\diag(\lambda_1,\lambda_1^{-1},-\lambda_1,-\lambda_1^{-1},\ldots,\lambda_r,\lambda_r^{-1},-\lambda_r,-\lambda_r^{-1},\mu_1,\mu_1^{-1},\ldots,\mu_s,\mu_s^{-1})
$ in $Sp(2n,\bar k)$ where $\mu_i^2=\pm 1$.
Such an element $t$ can be conjugated to $-t^{-1}$ by $s=\diag(\underbrace{M,\ldots,M}_r,\underbrace{N,\ldots,N}_s )\in Sp(2n,\bar k)$ where 
$$M=\left (\begin{array}{cccc}0&0& 0&-1 \\ 0&0& 1&0 \\ 0&-1&0&0\\ 1&0&0&0 \end{array}  \right)$$ 
and $s^2=-1$.
This concludes the proof of the theorem over $\bar k$.

We now complete the proof over $k$.
Let $t\in Sp(V)$, where  $V$ is a $2n$-dimensional vector space over $k$.
We first assume $t$ is real in $Sp(V)$.

First note that if $t_1\in Sp(V_1)$ and $t_2\in Sp(V_2)$, where $V_1$ and $V_2$ are vector space over $k$ of dimension $2n_1$ and $2n_2$ respectively, and if there exist  $g_1\in Sp(V_1)$ and $g_2\in Sp(V_2)$ such that $g_it_ig_i^{-1}=t_i^{-1}$ and $g_i^2=-1$, then $t_1\oplus t_2$ is conjugate to its inverse $t_1^{-1}\oplus t_2^{-1}$ by $g=g_1\oplus g_2$ in $Sp(V_1\oplus V_2)$ and $g^2=-1$.

Now let $t\in Sp(V)$ be real.
We write $\bar V$ for $V\tensor \bar k$ and $\bar V_{\alpha}=\{x\in \bar V \mid t(x)=\alpha x\}$, where $\alpha \in \bar k^{*}$.
Both $\bar V_1$ and $\bar V_{-1}$ are defined over $k$.
Let the subspaces $V_1$ and $V_{-1}$ of $V$ be the descents of $\bar V_1$ and $\bar V_{-1}$ respectively.
We note that the dimension of $V_{-1}$ is even, since the determinant of $t$ is $1$.
We now assume $\alpha\neq \pm 1$.
Let $\bar W_{\alpha},~W_{\alpha}$ and $k_{\alpha}$ be defined exactly as in the proof of Proposition~\ref{theoremkn2}. Then $\bar W_{\alpha}$ is a nondegenerate subspace of $\bar V$.
The subspace $W_{\alpha}$ is a direct sum of $m_{\alpha}$ two-dimensional subspaces over $k_{\alpha}$, which are stable under $t$ and 
$t$ restricted to each of these $2$-dimensional subspace is conjugate to $\diag\{\alpha,\alpha^{-1}\}$. 

By Lemma~\ref{sp2kbar}, there exists $g_{\alpha}\in Sp(W_{\alpha})$ with $g_{\alpha}^2=-1$ such that $g_{\alpha}t|_{W_{\alpha}}g_{\alpha}^{-1}=t|_{W_{\alpha}}^{-1}$.
The subspace $W_{\Gamma\alpha}=\oplus_{\sigma\in \Gamma} W_{\sigma\alpha}$ is defined over $k$ 
and the restriction of $t$ to this subspace is $t_{\Gamma\alpha}=\oplus_{\sigma\in\Gamma} t_{\sigma\alpha}$, where $t_{\sigma\alpha}=
t|_{W_{\sigma\alpha}}$.
Also $g_{\Gamma\alpha}=\oplus g_{\sigma\alpha}$ is defined over $k$ and conjugates $t$ to $t^{-1}$ on the subspace $W_{\Gamma\alpha}$.
We note that the $g_{\Gamma\alpha}^2=-1$.

Now we write $V=V_1\oplus V_{-1}\oplus_{\alpha\in\bar k^{*}}W_{\Gamma\alpha}$. 
Since the dimension of $V_{-1}$ is even, we may take $g_{-1}$  as the direct sum of $N=\left (\begin{array}{cc} 0&-1 \\  1&0 \\ \end{array}  \right)$ on this subspace, $\frac{1}{2}\dim(V_{-1})$ times.
Since $\dim(V)$ is even, it follows that dimension of $V_1$ is even and we may take $g_{1}$ as the direct sum of $N$,  $\frac{1}{2}\dim(V_{1})$ times, on this subspace.
Finally we take $g=g_1\oplus g_{-1}\oplus_{\alpha\in \bar k^{*}} g_{\Gamma\alpha}\in Sp(2n,k)$.
We have $g^2=-1$ and $gtg^{-1}=t^{-1}$. 

Now let us assume that $t$ is conjugate to $-t^{-1}$.
We follow the same proof as above except that we consider $\bar W_{\alpha}=\bar V_{\alpha}\oplus \bar V_{\alpha^{-1}}\oplus \bar V_{-\alpha}\oplus \bar V_{-\alpha^{-1}} $ when $\alpha^2\neq \pm 1$.
We construct $g_{\Gamma\alpha}$ using Lemma~\ref{sp4kbar} in this case.
The rest of the proof is along similar lines as above.    \qed
\vskip2mm

\noindent{\bf Remark :} We give an example to show that there are semisimple real elements in $Sp(4,k)$ which are not a product of two involutions.
Let $$J=\diag\left(\left (\begin{array}{cc} 0&-1 \\  1&0 \\ \end{array}  \right),\left (\begin{array}{cc} 0&-1 \\  1&0 \\ \end{array}  \right) \right)$$
be the matrix of the skew-symmetric (symplectic) form.
Then $Sp(4,k)=\{A\in GL(4,k)\mid \tr AJA=J\}$.
Let $S=\diag(\lambda,\lambda^{-1},\mu,\mu^{-1})\in Sp(4,k)$ with all diagonal entries distinct.
Then any element $T\in Sp(4,k)$, such that $TST^{-1}=S^{-1}$, is of the following type:
$$T=\diag\left(\left (\begin{array}{cc} 0&-a \\  a^{-1}&0 \\ \end{array}  \right),\left (\begin{array}{cc} 0&-b \\  b^{-1}&0 \\ \end{array}  \right) \right)$$ 
such that $T^2=-1$.
Hence $A$ is real semisimple but not a product of two involutions.

\subsection{Unitary Groups}
In this section we deal with unitary groups.
Let $K$ be a quadratic field extension of $k$.
Let $V$ be an $n$-dimensional vector space with a nondegenerate hermitian form $h$.
Then 
$$U(V,h)=\{t\in GL(V)\mid h(t(v),t(w))=h(v,w)\ \forall v,w\in V\}$$
is a $k$-group.
Let $\bar k$ be an algebraic closure of $k$.
We denote $\bar V=V\tensor_k \bar k$, a module over $K\tensor_k \bar k$.
We define $\bar h$ on $\bar V$ by base change of $h$ to $\bar k$.
Then $U(\bar V,\bar h)$ is an algebraic group defined over $k$ and $U(V,h)$ is the group of $k$-points of $U(\bar V,\bar h)$.
Let $\{e_1,\ldots,e_n\}$ be an orthogonal basis of $V$ with respect to $h$.
Let $h(e_i,e_i)=\alpha_i\in k$ and let $H=\diag(\alpha_1,\ldots,\alpha_n)$.
Then $U(V,h)\cong U(H)=\{A\in GL_n(K)\mid \tr AH\bar A=H\}$. We begin with a lemma for $V$ with $\dim(V)=2$.
\begin{lemma}
Let $V$ be a two dimensional vector space over $K$ with a nondegenerate hermitian form $h$.
Let ${e_1,e_2}$ be an orthogonal basis of $V$ with $h(e_i,e_i)=h_i$ and $H=\left (\begin{array}{cc} h_1&0 \\  0&h_2 \\ \end{array}  \right) $.
Let $A$ be any diagonal matrix in $U(H)$.
Then $A$ is real in $U(H)$ if and only if $h_1h_2\in N_{K/k}(K^*)$ and, in that case, it is strongly real.
\end{lemma}
\noindent {\bf Proof :}
Let $A=\left (\begin{array}{cc}\xi&0 \\ 0&\bar\xi \end{array}  \right)\in U(H)$.
Let $T$ be an element such that $TAT^{-1}=A^{-1}$.
Then $T$ is of the form: $T=\left (\begin{array}{cc}0&b \\ c&0\end{array}  \right)$ where $h_1b\bar b=h_2$ and $h_2c\bar c=h_1$.
Hence $A$ is real in $U(H)$ if and only if $h_1h_2\in N_{K/k}(K^*)$.
And, if the condition holds, we can take $T=\left (\begin{array}{cc}0&b \\ b^{-1}&0\end{array}  \right)$.  
This proves the result. \qed

\begin{theorem}\label{resultsinunitary}
Let $(V,h)$ be a hermitian space over $K$.
Let $t\in U(V,h)$ be a semisimple element.
Then, $t$ is real in $U(V,h)$ if and only if $t$ is strongly real.
\end{theorem}
\noindent{\bf Proof : }  
Let $t\in U(V,h)$ be a real semisimple element.
Let $g\in U(V,h)$ be such that $gtg^{-1}=t^{-1}$.
We base change to $\bar k$ and argue.
Since $t$ is real semisimple, we have a decomposition of $\bar V$ as follows:
$$\bar V=\bar V_1\oplus \bar V_{-1}\bigoplus_{\lambda\in \bar k^*}(\bar V_{\lambda}\oplus \bar V_{{\lambda}^{-1}})$$
where $\bar V_1, \bar V_{-1}$, $\bar V_{\lambda}$ and $\bar V_{{\lambda}^{-1}}$ are eigenspaces corresponding to eigenvalues $1, -1, \lambda$ and $\lambda^{-1}$ respectively.
Moreover, this decomposition is an orthogonal decomposition.
We denote the subspace $\bar V_{\lambda}\oplus \bar V_{{\lambda}^{-1}}$ by $\bar W_{\lambda}$.
It is easy to see that the conjugating element $g$ leaves $\bar W_{\lambda}$ invariant.
Since $\bar V_{\lambda}$ is nondegenerate, we can choose an orthogonal basis $\{e_1,\ldots,e_r\}$ for $\bar V_{\lambda}$.
We decompose $\bar W_{\lambda}$ in $t$ invariant planes as follows.
Let $P_i$ be the subspace generated by $\{e_i,g(e_i)\}$.
Then $\bar V_{\lambda}=P_1\oplus\ldots\oplus P_r$ is an orthogonal decomposition.
Moreover, $t$ leaves each of the $P_i$ invariant.
The element $n_i$ which maps $e_i$ to $g(e_i)$ and $g(e_i)$ to $e_i$, is a unitary involution conjugating $t|_{P_i}$ to its inverse. 
The element $\bar s=n_1\oplus\ldots\oplus n_r$ conjugates $t|_{\bar W_{\lambda}}$ to its inverse and is a unitary involution.

Let $W_{\lambda}$ be the sum of all Galois conjugates of $\bar W_{\lambda}$ and $s$ be the sum of all Galois conjugates of $\bar s$.
Then $W_{\lambda}$ is defined over $k$ and $t|_{W_{\lambda}}$ is conjugate to its inverse by the involution $s$ defined over $k$. 
This gives the decomposition of $V$ as $V=V_1\oplus V_{-1}\oplus_{\lambda} W_{\lambda}$ and we have proved that $t$ is a product of two involutions on each component.
Hence $t$ is strongly real.
\qed

\begin{corollary}
Let $t\in SU(V,h)$ be semisimple.
Suppose $n\not\equiv 2 \imod 4$.
Then $t$ is real in $SU(V,h)$ if and only if it is strongly real.
\end{corollary}
\noindent{\bf Proof:} The result follows by keeping track of determinant of the conjugating element in the proof of Theorem~\ref{resultsinunitary}.   \qed

{\bf Remarks :}
{\bf 1.} Let $K$ be a quadratic extension of $k$.
Let V be a two dimensional vector space over a field $K$ with a nondegenerate hermitian form $h$ defined as follows.
Let $\{e_1,e_2\}$ be a basis of $V$ such that $h(e_1,e_1)=1, h(e_2,e_2)=-1$ and $h(e_1,e_2)=0$.
In the matrix notation, the matrix of the form is $H=\left (\begin{array}{cc}1&0 \\ 0&-1 \end{array}  \right)$ and $U(H)=\{X\in GL_2(K)\mid \tr XH\bar X=H\}$.
Let $A=\left (\begin{array}{cc}\xi&0 \\ 0&\bar\xi \end{array}  \right)\in SU(H)$ where $\xi\neq \bar \xi$.
Then $A$ is semisimple.
Let $T\in GL_2(K)$ such that $TAT^{-1}=A^{-1}$.
Then $T$ is of the form $T=\left (\begin{array}{cc}0&b \\ c&0 \end{array}  \right)$.
Note that $A$ is real in $U(H)$ if and only if there exists $T=\left (\begin{array}{cc}0&b \\ c&0 \end{array}  \right)$ with $b\bar b=-1$ and $c\bar c=-1$.
The element $A$ is not strongly real in $SU(H)$.
For $T$ to be in $SU(H)$ we need $bc=-1$ and this implies $T^2=-1$.
Hence no involution conjugates $A$ to its inverse.
But if $K$ has an element $b$ such that $b\bar b=-1$, then $A$ can be conjugated to $A^{-1}$ by $T$ such that $T^2=-1$.
For example one can take $K=\mathbb Q(\sqrt 2)$ and $k=\mathbb Q$.

{\bf 2.}  Let V be a two dimensional vector space over $K$ with a hermitian form $h$ on it.
Let $K=k(\gamma)$.
Let $\{e_1,e_2\}$ be a basis of $V$ such that $h(e_1,e_1)=0, h(e_2,e_2)=0$ and $h(e_1,e_2)=\gamma=-h(e_2,e_1)$.
In the matrix notation, the matrix of the form is $H=\left (\begin{array}{cc}0&\gamma \\ -\gamma&0 \end{array}  \right)$ and $U(H)=\{X\in GL_2(K)\mid \tr XH\bar X=H\}$.
Let $A=\left (\begin{array}{cc}1&1 \\ 0&1 \end{array}  \right)\in SU(H)$.
Then $A$ is a unipotent element.
Let $T\in GL_2(K)$ be such that $TAT^{-1}=A^{-1}$.
Then $T$ is of the form $T=\left (\begin{array}{cc}a&b \\ 0&-a \end{array}  \right)$.
Note that $A$ is real in $U(H)$ if and only if there exists $T=\left (\begin{array}{cc}a&b \\ 0&-a \end{array}  \right)$ with $a\bar a=-1$ and $a\bar b-\bar ab=0$.
Here $T^2=a^2I$.
The element $A$ is not strongly real in $SU(H)$.
For if so, we would have $a^2=1$ and $a\bar a=-1$, which would imply that $\gamma$ is a square in $k$.
Hence no involution conjugates $A$ to its inverse.
But if $k$ has an element $a$ such that $a^2=-1$, then $A$ is conjugate to its inverse by $T$ such that $T^2=-1$.
For example one can take $k=\mathbb Q(\sqrt -1)$ and $K=\mathbb Q(\sqrt -1, \sqrt 5)$.

\appendix{}
\section{$G_2$ Revisited}

We take this opportunity to improve our result in~\cite{st} for all elements in $G_2$.
Let $G$ be a group of type $G_2$ defined over $k$.
In \cite{st}, we proved that a semisimple element in $G(k)$ is $k$-real if and only if it is strongly $k$-real and that unipotent elements in $G(k)$ are strongly $k$-real.
In this section we show that all real elements of $G(k)$ are strongly real in $G(k)$.
Since the proof is obtained by modifying the proof in the semisimple case, we shall refrain from repeating proofs of statements which 
are already there and provide appropriate references.
We follow the notation introduced in \cite{st}, Section 6.

\subsection{Reality in Groups of type $G_2$}

Let $G$ be a group of type $G_2$ defined over a field $k$ (of characteristic $\neq 2$).
 Then, there exists an octonion algebra $\C$ over $k$ such that $G\cong \Aut(\C)$ (\cite{se}, Chapter III, Proposition 5, Corollary). 
Let $t_0$ be an element of $G(k)$. 
We will also denote the image of $t_0$ in $\Aut(\C)$ by $t_0$.  
We let $V_{t_0}=ker(t_0-1)^8$. 
Then $V_{t_0}$ is a composition subalgebra of $\C$ with norm as the restriction of the norm on $\C$ (\cite{w2}). 
Let $r_{t_0}=\dim(V_{t_0}\cap \C_0)$, where $\C_0$ denotes the subspace of elements of trace $0$ in $\C$. 
Then $r_{t_0}$ is $1,3$ or $7$. 
We note that if $r_{t_0}=7$, the characteristic polynomial of $t_0$ is $(X-1)^8$ and $t_0$ is unipotent.
We have (\cite{st}, Theorem 6.3), 
\begin{lemma}\label{unipotentg2real}
Let $t_0\in G(k)$ be a unipotent element.
In addition, we assume $char(k)\neq 3$.
Then $t_0$ is strongly real in $G(k)$.
\end{lemma}
Let $L\subset\C$ be a quadratic \'{e}tale subalgebra. Let 
$$G(\C/L)=\{\phi\in G\mid \phi(x)=x,~\forall x\in L\}.$$
Recall from \cite{st}, when $L$ is a quadratic extension of $k$, $G(\C/L)\cong SU(L^{\perp}, h)$, for a nondegenerate hermitian form $h$ on the $3$ dimensional $L$-vector subspace $L^{\perp}$ of $\C$. When $L$ is split, $G(\C/L)\cong SL(3)$.  
\begin{lemma}\label{fixsubalgebra}
Let $t_0\in G(k)$ be an element which is not unipotent. 
Then, either $t_0$ leaves a quaternion subalgebra invariant or fixes a quadratic \'{e}tale subalgebra $L$ of $\C$ pointwise. 
\end{lemma}
\noindent{\bf Proof :} 
Since $t_0$ is not unipotent, from the above discussion, we see that $r_{t_0}$ is $1$ or $3$.
In the case $r_{t_0}=1$,  $L=V_{t_0}$ is a two dimensional composition subalgebra and has the form $V_{t_0}=k.1\oplus (V_{t_0}\cap \C_0)$, an orthogonal direct sum.
 Let $L\cap \C_0 = k.\gamma$ with $N(\gamma)\neq 0$. 
Since $t_0$ leaves $\C_0$ and $V_{t_0}$ invariant, we have, $t_0(\gamma) = \gamma$ and hence $t_0(x)=x\ \forall x\in L$, 
so that $t_0\in G(\C/L)$. When $r_{t_0}$ is $3$, the subalgebra $V_{t_0}$ is a quaternion algebra, left invariant by $t_0$. 
\qed
\vskip2mm
If $t_0$ leaves a quaternion subalgebra invariant, $t_0$ is strongly real in $G(k)$. 
This follows from Theorem 4 in \cite{w2} (see also Theorem 6.1 in \cite{st}). 
We discuss the other cases here, i.e., the fixed points of $t_0$ form a quadratic \'{e}tale subalgebra $L$ of $\C$.
\begin{enumerate}
\item The fixed subalgebra $L$ is a quadratic field extension of $k$ and 
\item the fixed subalgebra is split, i.e., $L\cong k\cross k$.
\end{enumerate}

By the above discussion, in the first case, $t_0$ belongs to $G(\C/L)\cong SU(L^{\perp},h)$ (Proposition 3.1 in \cite{st}).
We write $\C = L\oplus V$, where $V=L^{\perp}$ is a $3$-dimensional $L$-vector space with hermitian form $h$ induced by the norm on $\C$. 
In the second case, $t_0$ belongs to $G(\C/L)\cong SL(3)$ (Proposition 3.2 in \cite{st}). 
We denote the image of $t_0$ by $A$ in both of these cases. 
The characteristic polynomial $\chi_A(X)$ and the minimal polynomial $m_A(X)$ of $A$ will be refered to over $L$, in the first case and over $k$, in the second case.
We analyze further the cases depending on the characteristic polynomial of $A$.
We mention a result of Neumann here (\cite{n}, Satz 6 and Satz 8).
\begin{proposition}\label{neumannstreal}
Let the notation be as above. Let $t_0\in G(\C/L)$.  
Assume that the characteristic polynomial of $A$ is reducible and the minimal polynomial of $A$ is not of the form $(X-\alpha)^3$. 
Then $t_0$ is strongly real. 
\end{proposition}
\noindent We have the following,
\begin{theorem}\label{maintheorem}
Let $G$ be a group of type $G_2$ over a field $k$ of characteristic not $2$.
Let $t_0\in G(k)$ be an element which is not unipotent. 
Then, $t_0$ is real\index{real element} in $G(k)$ if and only if $t_0$ is strongly real in $G(k)$. 
In addition, if $char(k)\neq 3$ then every unipotent element in $G(k)$ is strongly real in $G(k)$. 
\end{theorem} 
\noindent{\bf Proof :} 
The assertion about unipotents in $G(k)$ is Lemma~\ref{unipotentg2real}.
In view of Lemma~\ref{fixsubalgebra} and discussion following the lemma, we need to consider the case when $t_0\in SU(V,h)$ or $t_0\in SL(3)$.
In these cases, we consider the characteristic polynomial $\chi_A(X)$ and the minimal polynomial $m_A(X)$ of $A$.
We first assume that $\chi_A(X)\neq m_A(X)$.
Hence degree of $m_A(X)$ is at most $2$ and $\chi_A(X)$ is reducible.
Clearly the minimal polynomial is not of the form $(X-\alpha)^3$.
Then by Proposition~\ref{neumannstreal}, $t_0$ is strongly real.
We take up the case of $A$ with $\chi_A(X)=m_A(X)$ below.
\qed
\vskip3mm
The result follows from the following 
\begin{theorem}\label{maintheoremg2semisimple} 
Let $t_0$ be an element in $G(k)$ and suppose $t_0$ fixes exactly a quadratic \'{e}tale subalgebra $L$ of $\C$ pointwise. 
Let us denote the image of $t_0$ by $A$ in $SU(V,h)$ or in $SL(3)$ as the case may be. 
Also assume that the characteristic polynomial of $A$ over $L$ in the first case and over $k$ in the second, is equal to the minimal polynomial of $A$. 
Then $t_0$ is conjugate to $t_0^{-1}$ in $G(k)$ if and only if $t_0$ is strongly real in $G(k)$.
\end{theorem}
\noindent{\bf Proof :} 
We distinguish the cases of both these subgroups below and complete the proof in the next two subsections, 
see Theorem~\ref{mainsuvh} and Theorem~\ref{counterexample}.
\qed

\begin{corollary}
Let characteristic $k\neq 2, 3$.
Then, an element $t\in G(k)$ is real in $G(k)$ if and only if $t$ is strongly real in $G(k)$.
\end{corollary}  
\subsection{$SU(V,h)\subset G$}\label{division}

We continue with notation introduced in the last section.
We assume that $L$ is a quadratic field extension of $k$. 
Let $t_0$ be an element in $G(\C/L)$ with characteristic polynomial of the restriction to $V$, equal to its minimal polynomial over $L$, i.e., $\chi_A(X)=m_A(X)$. 
We then have $G(\C/L)\cong SU(V,h)$. 
\begin{lemma}\label{sufix} 
Let $t_0$ be an element in $G(\C/L)$ which does not have a nonzero fixed point outside $L$. 
Suppose that $\exists g\in G(k)$ such that $gt_0g^{-1}=t_0^{-1}$. 
Then $g(L)=L$.
\end{lemma}
\noindent{\bf Proof :} 
Suppose $g(L)\not\subset L$. Then, as in the proof of Lemma 6.2 in \cite{st}, there exists $x\in L\cap \C_0$, a nonzero element, 
such that $g(x)\not\in L$. Since $t_0(x)=x$, it follows that $t_0(g(x))=g(x)$. Hence $t_0$ fixes $g(x)\not\in L$, a contradiction.  
\qed
\vskip3mm
We fix the basis for $V$ over $L$ introduced in the Section 6.1 in \cite{st}. 
Let us denote the matrix of $h$  with respect to this basis by $H=\diag(\lambda_1,\lambda_2,\lambda_3)$ where $\lambda_i=h(f_i,f_i)\in k^*$. 
Then $SU(V,h)$ is isomorphic to $SU(H)=\{A\in SL(3,L) \mid \tr AH\bar A=H\}$, where $a\mapsto \bar a$ is the nontrivial $k$-automorphism of $L$ and $\bar A$ is the matrix obtained by applying this automorphism to the entries of $A$.
\begin{theorem}\label{sucong2}
Let the matrix of $t_0$ be $A\in SU(H)$. 
Suppose that $t_0$ does not have a nonzero fixed point outside $L$. 
Then $t_0$ is conjugate to $t_0^{-1}$ in  $G(k)$, if and only if $\bar A$ is conjugate to $A^{-1}$ in $SU(H)$.
\end{theorem}
\noindent{\bf Proof :}
Let $g\in G(k)$ be such that $gt_0g^{-1}=t_0^{-1}$. 
By Lemma~\ref{sufix}, we have $g(L)=L$. Recall that $G(\C,L)\cong G(\C/L)\rtimes N$, where $N=<\rho>$ and $\rho$ is an automorphism of $\C$ with $\rho^2=1$ and $\rho$ restricts to the nontrivial automorphism of $L$. Using similar arguments as in the proof of Theorem 6.5 in \cite{st}, we conclude that $\bar A$ is conjugate to $A^{-1}$ in $SU(H)$.  Conversely, let $B\bar AB^{-1}=A^{-1}$ for some $B\in SU(H)$. 
Let $g'\in G(\C/L)$ be the element corresponding to $B$. 
Then $g'\rho$ conjugates $t_0$ to $t_0^{-1}$.
\qed

\begin{theorem}\label{mainsuvh}
Let $t_0$ be an element in $G(\C/L)$ which does not have a fixed point outside $L$ and let $A$ denote the image of $t_0$ in $SU(H)$. 
Suppose the characteristic polynomial of $A$ is equal to its minimal polynomial over $L$. 
Then $t_0$ is conjugate to $t_0^{-1}$, if and only if $t_0$ is a product of two involutions in $G(k)$. 
\end{theorem}
\noindent{\bf Proof :} 
From Theorem \ref{sucong2} we have, $t_0$ is conjugate to $t_0^{-1}$, if and only if $\bar A$ is conjugate to $A^{-1}$ in $SU(H)$. 
From Lemma 6.5 in \cite{st}, $\bar A$ is conjugate to $A^{-1}$ in  $SU(H)$ if and only if $A=A_1A_2$ with $A_1,A_2\in SU(H) $ and $\bar A_1 A_1=I=\bar A_2 A_2$. 
Now, from Proposition 6.1 in \cite{st}, it follows that $t_0$ is a product of two involutions.
\qed

\subsection{$SL(3)\subset G$}\label{split}

We continue here with proof of the Theorem~\ref{maintheoremg2semisimple}.
Let us assume now that $L\cong k \times k$. 
We have seen in \cite{st}, Section 3 that $G(\C/L)\cong SL(3)$. 
Let $t_0$ be an element in $G(\C/L)$ and denote its image in $SL(3)$ by $A$. 
 We assume that the characteristic polynomial of $A\in SL(3)$ is equal to its minimal polynomial over $k$.
\begin{lemma}
Let $t_0$ be an element in $G(\C/L)$ which does not have a fixed point outside $L$.
Suppose that $\exists h\in G=\Aut(\C)$, such that $ht_0h^{-1} =t_0^{-1} $. 
Then $h(L) = L$. 
\end{lemma}
\noindent{\bf Proof :} The proof is similar to that of Lemma~\ref{sufix}.

\qed

From Theorem~\ref{realsln} it follows that if $t_0$ is conjugate to $t_0^{-1}$ in $G(\C/L)\cong SL(3)$ then $t_0$ is strongly real.
Hence we may assume that $A$ is not real in $SL(3)$.
\begin{theorem}\label{splitconjugate}
Let $A$ be the matrix of $t_0$ in $SL(3)$ and assume that $A$ is not real in $SL(3)$.
Then $t_0$ is conjugate to $t_0^{-1}$ in $G=\Aut(\C)$, if and only if $A$ is conjugate to $\tr A$ in $SL(3)$. 
\end{theorem}
\noindent{\bf Proof :} 
Let $h\in G$ be such that $ht_0h^{-1}=t_0^{-1}$. 
Then, by the lemma above, $h(L) = L$. 
We may assume that (\cite{st}, Section 2) 
$$
\C =  \left \{ \left (\begin{array}{cc} \alpha &v \\  w&\beta \\ \end{array}  \right) \mid 
\alpha ,\beta \in k  ; v,w \in k^3 \right\}  \ \rm{with}\  L = \left \{ \left (\begin{array}{cc} \alpha & 0 
\\  0 &\beta \\ \end{array}  \right) \mid \alpha ,\beta \in k   \right\}. 
$$ 
Recall that $G(\C,L)\cong G(\C/L)\rtimes H$, where $H=<\rho>$ and $\rho$ is the automorphism of $\C$ which flips the diagonal and the anti-diagonal entries of a given element of the split octonion algebra $\C$ and the action of $SL(3)\cong G(\C/L)$ is as follows 
(see \cite{st}, Section 3): for $A\in SL(3)$ and for 
$$X=\left (\begin{array}{cc} \alpha & v \\ w & \beta \end{array} \right)\in \C,~~
AX=\left (\begin{array}{cc} \alpha & Av \\ {\tr A}^{-1}w & \beta \end{array} \right).$$
Hence, by the above lemma, 
$h\in G(\C/L)\rtimes H$. 
Since $A$ is not real in $SL(3)$, $h\notin G(\C/L)$.
Hence $h = g\rho$ for some $g \in G(\C/L)$. 
Let $B$ denote the matrix of $g$ in $SL(3)$. 
Then, a computation same as in the proof of Theorem 6.7 of \cite{st}, shows
$$
ht_0 h^{-1} = t_0^{-1} \Leftrightarrow A =B \tr A B^{-1}.
$$ 
Therefore $t_0$ is conjugate to $t_0^{-1}$ in $G(k)$ if and only if $A$ is conjugate to $\tr A$ in $SL(3)$.
\qed

\begin{lemma}\label{consln}
Let $A$ be a matrix in $SL(n)$ with its characteristic polynomial equal to its minimal polynomial. 
Then $A$ is conjugate to $\tr A$ in $SL(n)$ if and only if $A$ is a product of two symmetric matrices in $SL(n)$.
\end{lemma}
\noindent{\bf Proof :} 
The proof is exactly same as the proof of Lemma 6.10 in \cite{st}.   
\qed
\begin{theorem}\label{counterexample}
Let $t_0\in G(\C/L)$. 
Assume that the characteristic polynomial of the matrix $A$ of $t_0$ in $SL(3)$ is equal to its minimal polynomial. 
Then, $t_0$ is conjugate to $t_0^{-1}$ in $G=Aut(\C)$ if and only if $t_0$ is a product of two involutions in $G(k)$.
\end{theorem}
\noindent{\bf Proof :} 
First, let $t_0$ be real in $G(\C/L)$. 
Then, $A$ is real in $SL(3)$ and hence it is strongly real (see Theorem~\ref{realsln}).
Thus the element $t_0$ is strongly real in $G(k)$.
Now we assume $t_0$ is not real in $G(\C/L)$, i.e., $A$ is not real in $SL(3)$.
In this case, the element $t_0$ can be conjugated to $t_0^{-1}$ in $G(k)$ if and only if, $A$ can be conjugated to $\tr A$ in $SL(3)$ (Theorem~\ref{splitconjugate}). 
This is if and only if, $A$ is a product of two symmetric matrices in $SL(3)$ (Lemma~\ref{consln}). 
By Proposition 6.5 in \cite{st}, this is if and only if $t_0$ is a product of two involutions in $\Aut(\C)$.
\qed

\vskip2mm
{\bf Acknowledgment :} The aunthors are indebted to the referee for his/her invaluable suggestions. 
We take this opportunity to thank Prof. T. A. Springer 
and Prof. Dipendra Prasad for their help and encouragement.



\begin{thebibliography}{99}

\bibitem[D71]{d} D. \v Z. Djokovi\'c, \emph{``The product of two involutions in the unitary group of a hermitian form''}, Indiana Univ. Math. J. 21 1971/1972, 449-456.
\bibitem[E79]{e} Erich W. Ellers, \emph{``Products of two involutory matrices over skewfields''}, Linear Algebra Appl. 26 (1979), 59-63. 
\bibitem[KMRT98]{kmrt} Max-Albert Knus, Alexander Merkurjev, Markus Rost and  Jean-Pierre Tignol, \emph{``The book of involutions''}, American Mathematical Society Colloquium Publications 44; American Mathematical Society, Providence, RI, 1998.
\bibitem[KN87]{kn} F. Kn\"uppel and K. Nielsen, \emph{``Products of involutions in ${\rm O}\sp +(V)$''},  Linear Algebra Appl.  94  (1987), 217-222.
 \bibitem[N90]{n} A. Neumann, \emph{``Bedingungen f\"{u}r die Zweispiegeligkeit der Automorphismengruppen von Cayleyalgebren''}, Geometriae Dedicata 34 (1990), no. 2, 145-159.
\bibitem[Pr98]{pr1} D. Prasad, \emph{``On the self-dual representations of finite groups of Lie type''},  J. Algebra  210  (1998),  no. 1, 298-310.
\bibitem[Pr99]{pr2} D. Prasad, \emph{``On the self-dual representations of a $p$-adic group''},  Internat. Math. Res. Notices   (1999),  no. 8, 443-452.
\bibitem[RS90]{rs} R. W. Richardson and T. A. Springer, \emph{``The Bruhat order on symmetric varieties''}, Geom. Dedicata  35 (1990), no. 1-3, 389-436.
\bibitem[S65]{s} R. Steinberg, \emph{``Regular elements of semisimple algebraic groups''},  Inst. Hautes Études Sci. Publ. Math. No. 25 1965, 49-80.
\bibitem[SS68]{ss} T. A. Springer, R. Steinberg, \emph{``Conjugacy classes''}, 1970  Seminar on Algebraic Groups and Related Finite Groups (The Institute for Advanced Study, Princeton, N.J., 1968/69),  pp. 167-266 Lecture Notes in Mathematics, Vol. 131 Springer, Berlin.
\bibitem[Se97]{se} J. P. Serre, \emph{``Galois cohomology''}, Springer-Verlag, Berlin, (1997).
\bibitem[ST05]{st} A. Singh, M. Thakur, \emph{``Reality Properties of Conjugacy Classes in $G_2$''}, Israel Journal of Mathematics 145 (2005), 157-192.
\bibitem[TaZ59]{taz} O. Taussky and H. Zassenhaus, \emph{``On the similarity transformation between a matrix and its transpose''}, Pacific Journal of Mathematics  9 (1959), 893-896.
\bibitem[TZ05]{tz} P. H. Tiep, A. E. Zalesski, \emph{``Real conjugacy classes in algebraic groups and finite groups of Lie type''}, J. Group Theory 8 (2005) no. 3, 291-315.
\bibitem[W66]{w1} M. J. Wonenburger, \emph{``Transformations which are products of two involutions''}, J. Math. Mech. 16 (1966), 327-338.
\bibitem[W69]{w2} M. J. Wonenburger, \emph{``Automorphisms of Cayley algebras''}, Journal of Algebra 12 (1969), 441-452.
\end{thebibliography}
\end{document}